\documentclass[11pt]{amsart}

\usepackage[english]{babel}

\usepackage{amsmath,amsfonts,amsthm,amssymb}
\usepackage{color}

\def\R{{\mathbb R}}

\def\P{{\bf P}}
\def\E{{\bf E}}
\def\Z{{\mathbb Z}}

\def\N{{\mathbb N}}

\def\Np{\mathcal N_{poss}}

\def\eps{{\varepsilon}}
\def\n{{\vec{n}}}

\newtheorem{defi}{Definition}
\newtheorem{theo}{Theorem}
\newtheorem{prop}{Proposition}
\newtheorem{cor}{Corollary}

\author[O. Adelman]{Omer ADELMAN}
\address{Laboratoire de Probabilit\'es et Mod\`eles Al\'eatoires, Universit\'e
Paris 6, 4 place Jussieu, 75252 Paris cedex 05}
\email{adelman@ccr.jussieu.fr, enriquez@ccr.jussieu.fr}

\author[N. Enriquez]{Nathana\"el ENRIQUEZ}

\title[What a single trajectory tells]{
Random walks in random environment:\\
What a single trajectory  tells}

\begin{document}
\maketitle {\bf Abstract}: {\small We present a procedure that determines the
law of a random walk in an iid random environment as a function of a
single ``typical"
trajectory. We indicate when the trajectory characterizes the law of the
environment, and we say how this law can be determined. We then show
how independent trajectories having the distribution of the original walk
can be generated
as functions of the single observed trajectory.}

\bigskip
\section{Introduction}

Suppose you are given a ``typical" trajectory of a random walk in
an iid random environment. Can you say what the law of the
environment is on the basis of the information supplied by this
single trajectory? Can you determine the law of the walk? Such
questions may arise if one intends to use the random environment
model in applications.

These questions are essentially pointless if the group is finite
(in which case the environment at each of the finitely many sites
that happen to be visited infinitely many times can of course be
determined, but it is hard to say much more). So we assume that
the group is infinite, and we go a little further: we assume that
the (random) set of sites visited by the walk is almost surely
infinite. (See remark 5.1.)

Questions of this kind have been studied in the context of random
walks in random scenery by Benjamini and Kesten \cite{BK}, L\" owe
and Matzinger \cite{LM}, and  Matzinger \cite{M}.

In the case of an iid random environment, the information furnished
by a single (``typical") trajectory tells us whether the walk is
recurrent; indeed, one can show that one of the events \{each visited site
is visited infinitely many times\}, \{no site is visited infinitely
many times\}
is an almost sure event. (Cf Kalikow \cite{K}.)

Now, if the walk is recurrent, the problem is quite simple: we can
know much more than the {\emph {law}} of the environment, because
we find the environment {\emph {itself}} at each visited site,
which is given by the frequency of each possible jump from this
site. In the transient case, the ``na\"{\i}ve" approach consisting
of doing statistics on sites which have been visited many times
cannot be utilized directly, since the assumption of being at a
site which has been frequently visited introduces a bias on the
environment at that site, which should encourage jumps to sites
from which it is be easier to come back (loosely speaking, close
sites).

We present a  procedure that eliminates any source of bias,
  collecting information on sites displaying some specified ``histories".
Each such ``history" which $can$ be encountered $is$ almost surely
encountered infinitely often (Proposition 3).
  This is combined with an interpretation of the
process as an edge-oriented reinforced random walk (cf
Enriquez and Sabot \cite{ES}), allowing us to find
the exact law of the process. Now, there may exist ``bad"
transitions: if the walk jumps from a site along a ``bad" transition, it will
never get back to that site. If the set of these ``bad" transitions is empty or
has just one element, we can find the distribution of the environment (Theorem
1).

Finally, we show how countably many independent trajectories
having the distribution of the original walk can be generated by
concatenating steps of the observed trajectory. The algorithm is
purely ``mechanical": it does not imply any computation, and, in
particular, the knowledge of the law of the walk (or of the
environment) is not needed.


\section{Framework and notations}

Our ``canonical" process $X:=(X_n)_{n\geq0}$ walks on a group $G$. We denote
by $(\mathcal F_n)_{n\geq0}$ the natural filtration of $X$.

    We assume that the group \,$G$\, is  Abelian, although this is
{\em{never}} used
in our arguments. Its only utility is in the possibility of
writing things like \,$x-y=e$\, or \,$x=y+e$\, indifferently.

We use the additive notation, and the identity element of \,$G$\,
is denoted by \,$0$.

We assume moreover that the group \,$G$\, is countable. This
assumption can be dispensed with$-$see remark 5.2$-$but we feel that it
renders the  reading easier and
the discussion more tractable. It does not affect the core of the
argument.

We denote by  \,$\N$\, the set \,$\{1,2,\dots\}$.

\subsection{Random walks in random environment\textcolor{white}{.}}~%

We denote by \,${\mathcal P}$\, the set of non-negative families
\,$p:= (p_e)_{e\in G}$\, such that \,$\sum_{e\in
G}p_e=\nolinebreak 1$. The environment at the site \,$x$,\,
$\nu(x):=(\nu(x,e))_{e\in G}$,\, is a random element of
\,${\mathcal P}$. We assume that the environments at sites are iid
\,${\mathcal P}$-valued random variables with common distribution
\,$\mu$.

We let \,$\nu_e:=\nu(0,e)$.

The random environment \,$\nu:=(\nu(x))_{x\in G}$\, is a random
element of \,${\mathcal P}^{G}$,\, and it is governed by the
probability measure \,$\mu^{\otimes G}$.\,

For all \,$\pi=((\pi(x,e))_{e\in G})_{x\in G}\in {\mathcal P}^{G}$,\, let
\,$Q_{\pi}$\, be the probability measure  under which
$X$ is a $G$-valued Markov chain started at 0 whose transition
probability from $x$ to $x+e$ is $\pi(x,e)$ \,($x, e\in G$).

The law of the random walk in
random environment (or the so-called ``\emph{annealed}" law) is the probability
measure \,$\P^\mu=\int Q_{\pi} \, \mu^{\otimes G}(d\pi) \ (=\E[Q_{\nu}]$,\
  $Q_{\nu}$ being what is usually called the ``quenched" law).

We recall our ``infinitude assumption" according to which the
(random) set of sites visited by the walk, \,$S:=\{X_n\,|\,n\ge
0\}$,\, is \,$\P^\mu$-almost surely infinite. (This implies, of
course, that the group \,$G$\, itself is infinite.)\medskip

We let $E$ denote the set of those $g\in G$ such that the
probability of the event $\{\nu_g>0\}$ is strictly positive. (It
is easy to see that the random set $\{X_{n+1}-X_n\,|\,n\ge 0\}$ is
$\P^\mu$-almost surely exactly $E$.) We then partition $E$ into two sets,
$R$ and $T$, defined as follows.

\noindent$\bullet$\ \  $R$ is the set of elements $r$ of $E$ that
can be written as $-(e_1+...+e_n)$ where $(e_i)_{1\leq i\leq n}$
is a finite nonempty sequence of elements of $E$. It is easy to
see that  $r\in R$ if and only if $\P^\mu(X_1=r \textrm { and, for some }
n>1,\ X_n=0)>0$\,; and Proposition 3 below implies that if $r\in
R$, then the random set $\{n\,|\,X_{n+1}=X_n+r {\textrm{  and, for
some }} k>0,\ X_{n+k}=X_n\}$ is $\P^\mu$-almost surely infinite. It is
therefore quite easy to identify $R$ when observing a single
trajectory.

\noindent$\bullet$\ \   $T$ is the complement of $R$ in $E$. It
represents the ``possible" transitions which do not allow a return
to the original site.\medskip

\subsection{Histories\textcolor{white}{.}}~%

We start with some definitions.

\begin{defi}
The \,{\rm  history  of the site $x$ at time $n$}, which we denote
by \,$H(n,x)$,\, is the random finite sequence of elements  of
\,$G$\, defined by the successive moves of the process from the
site \,$x$\, before time \,$n$. More formally, \,$H(0,x)$\, is the
empty sequence \,$(\,)$;\, $H(n+1,x)=H(n,x)$\, if \,$X_n\neq x$;\,
and, if \,$X_n= x$, \,then \,$H(n+1,x)$\, is the finite sequence
obtained by adjoining \,$X_{n+1}-X_n$\, as a new rightmost term to
\,$H(n,x)$.
\end{defi}

Let us denote by \,$(\Z_+^G)_0$\, the set of families
\,$(n_g)_{g\in G}\in \Z_+^G$ \,with a finite number of non null
terms.

\begin{defi}
The \,{\rm unordered history of the site \,$x$\, at time \,$n$},\,
denoted by \,$\vec{N}(n,x):=(N_g(n,x))_{g\in G}$,\, is a random
element of \,$(\Z_+^G)_0$\, where, for all \,$g\in G$,\,
$N_g(n,x)$\, is the random number of moves  from the site \,$x$\,
to \,$x+g$\, before time \,$n$. In other words,
\,$N_g(n,x)=\sum_{l=0}^{n-1}1_{\{X_l=x,\,X_{l+1}-X_l=g\}}$.

Also, the {\rm local unordered history at time $n$} is the unordered history of
the site $X_n$ at time $n$, $\vec{N}(n):= \vec{N}(n,X_n)$.
\end{defi}

\subsection{Reinforced random walks\textcolor{white}{.}}~%

An edge-oriented reinforced random walk consists of a discrete
random process whose transition probabilities are functions of the
number of each type of move in the history of the process that has
been made from the site currently occupied. A good point of view,
in order to get a non-biased procedure of reconstitution of the
environment, is to view the random walk in random environment as an
edge-oriented reinforced random walk. It is the essence of the
easy part of the result of \cite{ES} (the other part examines the
conditions on a reinforced random walk to correspond to a RWRE).

We introduce the reinforced random walks by the following
definitions.
\begin{defi}
A \,{\rm reinforcement function} is a function
$$V:(\Z_+^G)_0\to {\mathcal P}$$
$$\n=(n_g)_{g\in G}\mapsto V(\n):=(V_e(\n))_{e\in G}$$
\end{defi}

\begin{defi}
We call \, {\rm edge-oriented reinforced random walk with reinforcement
function} $V$ the random walk defined by the law $\P^V$ on the
trajectories starting at $0$ given by
$$\P^V(X_{n+1}-X_n=e\,|\,{\mathcal F}_n)=V_e(\vec{N}(n)).$$

\end{defi}


\section{Tools}

\subsection{RWRE as an edge-oriented reinforced random walk%
\textcolor{white}{.}}~%

We can now state the result of  Enriquez and Sabot \cite{ES}:
\begin{prop}
        The annealed law \,$\P^\mu$\, of the RWRE  coincides
with the law \,$\P^V$\, of the reinforced random walk  whose
reinforcement function \,$V$\, satisfies, for all \,$e\in G$,
$$ V_e(\n)=
{\E[\nu_e\prod_{g\in G}\nu_g^{n_g}]\over \E[\prod_{g\in
G}\nu_g^{n_g}]}
$$
whenever \,$\n\in(\Z_+^G)_0$\, such that \,$\E[\prod_{g\in
G}\nu_g^{n_g}]>0$.
\end{prop}

In order to be self-contained we recall the proof of this
proposition.

\noindent{\sl{P\tiny{ROOF}.\ }} %
For any  $x,e$ in $G$, for all $n\in\N,$ $\P^\mu$-almost
everywhere on the event $\{X_n=x\}$,

$$\P^\mu(X_{n+1}=x+e\,|\,{\mathcal F}_n)
={\E[\nu(x,e)\prod_{y\in G}\prod_{g\in G}\nu(y,g)^{n_g(n,y)}]\over
\E[\prod_{y\in G} \prod_{g\in G}\nu(y,g)^{n_g(n,y)}]}$$

Now using the independence of the random variables $\nu(y)$ for
different sites $y$, the terms depending on $\nu(y)$ for $y\neq x$
cancel in the previous ratio,  and we get the result. \qed

The following result is an analogue of the strong Markov property
for reinforced random walks.

\begin{prop}
Let $X$ be a \, {\rm reinforced random walk with reinforcement
function} $V$, and let $T$ be a stopping time with respect to the
natural filtration of $X$. Assume $T$ is  almost surely finite.
Then
$$\P^V(X_{T+1}-X_T=e\,|\, {\mathcal F}_T)=V_e(\vec{N}(T))\qquad
\P^V-a.s.$$
\end{prop}
The proof is obtained in an obvious way, by considering the events
$\{T=n\}$.

\subsection{A zero-one result\textcolor{white}{.}}~%

The following zero-one result happens to be quite useful.

\begin{prop}
Let $(r_1,...,r_l)$ a finite (eventually empty) sequence of
elements of $R$. Let $S_{(r_1,...,r_l)}$ be the random set $\{x\in
G\,|\,\exists n\geq 0,\, H(n,x)=(r_1,...,r_l)\}$. Then
$S_{(r_1,...,r_l)}$ is either $\P^\mu$-almost surely empty or $\P^\mu$-almost
surely infinite.
\end{prop}
\noindent{\sl{P\tiny{ROOF}.\ }} %
Suppose that $S_{(r_1,...,r_l)}$ is not $\P^\mu$-almost surely empty.

This implies that there exists a list of transitions
$$
L:=(r_{1},e_{1,1},e_{1,2},...,e_{1,k_1},r_{2},
e_{2,1}...,e_{2,k_2}, r_{3},.....,e_{l-1,k_{l-1}},r_{l})$$ such
that, for all  $m\in \{1,\dots, l-1\}$,
    $$r_m+\sum_{i=1}^{k_m}e_{m,i}=0\qquad \mbox{and, for all} \,\, k\in
\{1,\dots, k_m-1\},\quad r_m+\sum_{i=1}^{k}e_{m,i}\neq0$$ and
$$\gamma:=\E[\prod_{k=1}^l\nu_{r_k}]
\E[\prod_{m=1}^{l-1}\prod_{i=1}^{k_m-1}\nu(r_m+e_{m,1}+...+e_{m,i}\,
,\,e_{m,i+1})]>0.$$ (Note that if $r_m=0$, then $k_m=0$.)

Let $q:=l+k_1+...+k_{l-1}$ be the length of the list $\,L$.

For all $k\in\{1,\dots, q\}$, denote by $y_k$ the $k$-th  term of the
list $L$. For all $k\in\{1,\dots, q+1\}$, set $x_k:=\sum_{i=1}^{k-1}y_i$\,.
($\,x_1= 0$.)

Now consider the list $(g_1,g_2,...)$ of  newly visited sites in
their order of appearance. So ${\mathcal S}= \{g_1,g_2,...\}$. By the
assumption made in the introduction, ${\mathcal S}$ is almost
surely infinite.

To any integer $n\geq 1$, we associate a random integer $i(n)$
defined by
$$i(n):=\min\{i\geq1\,|\, \exists k\in\{1,\dots,q\},\,\, g_i=g_n+x_k\}. $$
We denote by $k(n)$ the random smallest integer $m\ge1$  such that
$g_{i(n)}=g_n+x_{m}$. (Obviously, $k(n)\le q$.) The sequence
$(g_{i(n)})_{n\geq1}$ takes infinitely many values (since the
infinite set ${\mathcal S}$ is included in $\{g_{i(1)},
g_{i(2)},...\}-\{x_1,...,x_q\}$). Now, for any $i\geq1$, we denote
by $T_i$ the hitting time of $g_i$ by the walk. By definition of
$i(n)$, none of the sites $g_{i(n)}-x_{k(n)}+x_j\quad(1\leq j\leq
q)$ is visited by the trajectory before time  $T_{i(n)}$.

As a result,  there exist  infinitely many sites $g'_1,g'_2,\dots$
visited by the trajectory (enumerated in their order of
appearance) such that, for some $k\in\{1,\dots, q\}$, if $T'_n$
denotes the hitting time of $g'_n$, then none of the sites
$g'_{n}-x_{k}+x_j\quad(1\leq j\leq q)$ is visited by the
trajectory before time $T'_n$. We denote by $k'(n)$ the least such
integer $k$. $T'_n$ are clearly stopping times.

For all $n\geq1$, let $\psi_n$ be the Bernoulli variable that
equals 1 if and only if
$$X_{T'_n+i}=
\left\{
{\begin{array}{ll}
g'_{n}-x_{k'(n)}+x_{i+k'(n)}&\qquad if \quad i\in\{1,\dots ,q-k'(n)\},\\
g'_n-x_{k'(n)}+x_{i-k'(n)-q}&\qquad  if \quad i\in\{q+1-k'(n),\dots
,2q-k'(n)+1\}.\end{array}}\right.
$$

(Otherwise, $\psi_n$
equals 0.)

Observe that for all $n$, if $\psi_n=1$, then
$X_{T'_n}-x_{k'(n)}\in S_{(r_1,...,r_l)}$.

Due to the fact that the prescribed path the process has to follow
during the time period \,$[T'_n, T'_n+2q-k'(n)+1]$\, in order to
satisfy \,$\psi_n=1$\, is a path that does not intersect the
trajectory before \,
$T'_n$,
$$ \P(\psi_n=1\,|\,{\mathcal F}_{T'_n}
)\geq\E[\prod_{k=1}^{q}(\nu(x_k,y_k))^2]\geq
\E[\prod_{k=1}^{q}\nu(x_k,y_k)]^2 = \gamma^2.$$
Let
$\xi_{n}:=\psi_{(2q+1)n}$ and $\tau_n=T'_{(2q+1)n}$ \ ($n\geq1$).

For all $n$, $\xi_1,\dots,\xi_n$ are measurable with respect to
$\mathcal F_{\tau_{n+1}}$.

It is now obvious that for all $n,k\geq1$
$$\P(\xi_{n+1}=\dots=\xi_{n+k}=0)=$$
$$\P(\xi_{n+1}=0)\cdot\P(\xi_{n+2}=0\,|\,
\xi_{n+1}=0)\cdots\P(\xi_{n+k}=0\,|\,
\xi_{n+1}=\cdots=\xi_{n+k-1}=0)\leq(1-\gamma^2)^k.$$ Therefore, %
almost surely, infinitely many of the  $\psi_n$'s are equal to 1.
This implies that $S_{(r_1,...,r_l)}$ is almost surely infinite.
\qed

\noindent{\sl{R\tiny{EMARK}.\ }} If $G$ equals $\Z^d$, the notion
of convexity can be exploited in a proof slightly different from the one given
above.

We deduce that the sets $R$ and $T$ can be ``viewed" on the trajectory:
\begin{cor} The sets  $R$ and $T$ are such that

\noindent
$R\,{\stackrel{a.s}{=}}\,\{g\in E\,|\, S_{(g)} \mbox{ is
infinite}\}\,{\stackrel{a.s}{=}}\,\{g\in E\,|\, S_{(g)} \neq\emptyset\}$
\ and
\  $T\,{\stackrel{a.s}{=}}\,\{g\in E\,|\, S_{(g)} =\emptyset\}$.
\end{cor}

Let $S_\n$ denote the random set \,$\{x\in
G\,|\, \exists n\geq 0,\, \vec N(n,x)=\n\}$. An easy corollary of the above
proposition is the following analogous result concerning unordered
histories.\medskip

\noindent{\bf{Proposition 3'.}}{\it{%
\ \ If\, $\n\in(\Z_+^G)_0$, then\, $S_\n$\,
   is either $\P^\mu$-almost surely empty or $\P^\mu$-almost surely
infinite.}}\bigskip

We now distinguish a particular subset of $(\Z_+^G)_0$,
$$\mathcal N_{poss}:=\{\n\in(\Z_+^G)_0\,|\,  S_\n\neq\emptyset\quad \P^\mu-
a.s.\}.$$
Loosely speaking, \,$\mathcal N_{poss}$\, is the set of
``possible" unordered histories for sites that are presently occupied.
An element $\n=(n_g)_{g\in G}$ of $(\Z_+^G)_0$ belongs to
$\mathcal N_{poss}$ if $n_g=0$ whenever $g\notin R$ and if, moreover, it
satisfies (any one of) the three following equivalent conditions:

\noindent(a)\quad $\P^\mu(\exists n\geq0\,|\,\vec N(n,0)=\n)>0$;

\noindent(b)\quad the random set $S_\n$ is almost surely infinite;

\noindent(c)\quad $\E[\prod_{r\in R}\nu_r^{n_r}]>0$.

Note that, by (c) and Proposition 1, the law of the process is determined by
the restriction of the reinforcement function to the set $\mathcal N_{poss}$.

\section{What a single trajectory  tells}

In the sequel we assume that the law of the process $X$ is $\P^\mu$ (or,
equivalently, $\P^V$).

\subsection{Determining the law of the walk\textcolor{white}{.}}~%

   As noted in the introduction, ``straightforward"
computation based on the frequencies
of transitions from sites visited ``many" times is not reliable.
What we do here, instead, is collecting information on sites displaying
some specified histories (or  specified \emph{unordered} histories).

\noindent For any $\n\in\mathcal N_{poss}$, let
$T_i^{\n}\quad(i\geq 1)$
   be the successive times where the unordered history of the
   currently occupied site is
$\n$:

\noindent$\bullet$\ \ $T_0^{\n}:= \inf\{k\geq0\,|\,
\vec{N}(k)=\n\}$

\noindent$\bullet$\ \ $\forall i\geq0, \, T_{i+1}^{\n}:= \inf\{k>
T_i^{\n} \,|\, \vec{N}(k)=\n\}$\\\medskip%
(We ignore the case (happening on the negligible event $\{S_\n$ is
finite\}) where some $T_i^{\n}$ is infinite.)

\begin{prop}
   The $E$-valued random variables $\Delta_i^{\n}:=X_{T_i^{\n}+1}-X_{T_i^{\n}}$
($i\geq1$, $\n\in\Np$) are independent. Also, for all $\n\in\mathcal N_{poss}$,
the random variables $\Delta_i^{\n}$ have the same law:
$$\forall e\in E,\quad\P(\Delta_i^{\n}=e)=V_e(\n).$$
\end{prop}

\noindent{\sl{P\tiny{ROOF}.\ }} %
Let $\Theta_i^{\n}$ ($i\geq1$, $\n\in(\Z_+^G)_0$) be
independent random variables such that for all $e\in E$,
$$\P(\Theta_i^{\n}=e)=V_e(\n)   $$
Now we consider the process $(Y_n)_{n\geq0}$ :
$$ Y_0=0,\quad  Y_{n+1}=Y_n+\Theta_{\tau(n)}^{\vec{N}_y(n)},$$
where $ \vec{N}_y$ is to the process $Y$ what $\vec{N}$ is to the
process $X$, and $$\tau(n):=
card\{j\in\{0,\cdots,n\}\,|\,\vec{N}_y(j)=\vec{N}_y(n)\}.$$
$$\begin{array}{rl}%
\P(Y_{n+1}=Y_n+e\,|\,\sigma(Y_{k},  k\leq n))\!\!\!&=
\P(\Theta_{\tau(n)}^{\vec{N}_y(n)}=e\,|\,\sigma(Y_{k},  k\leq
n))\vspace{0.3cm}\\%
&=\E[\sum_{l\geq0,\vec{m}\in(\Z_+^G)_0}1_{\tau(n)=l,\vec{N}_y(n)=\vec{m}}
1_{\Theta_{l}^{\vec{m}}=e}\,|\,\sigma(Y_{k},  k\leq
n)]
\end{array}$$
But on the event
$A_{l,\vec{m},n}:=\{\tau(n)=l,\vec{N}_y(n)=\vec{m}\}$,
$\Theta_{l}^{\vec{m}}$ is independent of $\sigma(Y_{k},  k\leq
n)$. Thus,
$$\begin{array}{rl}%
\P(Y_{n+1}=Y_n+e\,|\,\sigma(Y_{k},  k\leq n))\!\!\!&=
\sum_{l\geq0,\vec{m}\in(\Z_+^G)_0}1_{\tau(n)=l,\vec{N}_y(n)=\vec{m}}
\E[1_{\Theta_{l}^{\vec{m}}=e}]\vspace{0.3cm}\\
&=\sum_{l\geq0,\vec{m}\in(\Z_+^G)_0}1_{\tau(n)=l,\vec{N}_y(n)=\vec{m}}
V_e(\vec{m})\vspace{0.3cm}\\
&=V_e(\vec{N}_y(n))
\end{array}$$
Consequently, the  two processes $X$ and  $Y$ have the law. But
$\Delta_i^{\vec{n}}$ is to the process $(X_n)_{n\geq0}$  exactly
what $\Theta_i^{\vec{n}}$ is to the process $Y$. The result
follows. \qed

We deduce from this proposition the following corollary, which
describes the construction of the reinforcement function on
$\mathcal{N}_{poss}$ or, equivalently (by Proposition 1), of the
annealed law:

\begin{cor}
If $\n\in\mathcal{N}_{poss}$, then almost surely, for all $ e\in
E$, \,
$${1_{\Delta_1^{\n}=e}+...+1_{\Delta_m^{\n}=e}\over
m}\to V_e(\n)\quad as\quad{m\to\infty}.$$
\end{cor}

\subsection{The law of the environment\textcolor{white}{.}}~%

The next result follows easily.

\begin{theo}
(a) A single trajectory determines almost surely the moments
of the form
$\E_\mu[\nu_{r_1}^{n_1}\cdots\nu_{r_k}^{n_k}\nu_t^\eps]$ for all
$r_1,..., r_k\in R$, $t\in T$, $n_1,..., n_k\in \Z_+$, $\eps=0$ or
$1 $. Moreover, if these moments coincide for two distinct
environment distributions,  the induced RWRE have the same
annealed law (and, consequently, two such environment
distributions cannot be distinguished).

(b) If  \,{\rm card} $T= 0$ or $1$, a single trajectory determines almost
surely the distribution of the environment.

\end{theo}

\noindent{\sl{P\tiny{ROOF}.\ }} %
(a) By Corollary 2, the restriction of $V$ to $\mathcal{N}_{poss}$ is almost
surely determined by a single trajectory. So almost surely, for all
$\n=(n_g)_{g\in G}\in \mathcal{N}_{poss}$
and for all $e\in G$, a single trajectory determines the moments
$\E[(\prod_{g\in G}\nu_{g}^{n_g})\cdot\nu_e]$. Moreover, all the
other moments of
the type $\E_\mu[\nu_{r_1}^{n_1}\cdots\nu_{r_k}^{n_k}\nu_t^\eps]$ ($r_j\in R$,
$t\in T$) are zero.
Finally, the restriction of $V$ to $\mathcal{N}_{poss}$
determines the law of the process.

(b) If card\,$T= 0$, we get all the moments of the $\nu_e$'s. Since these
variables are compactly supported variables,  this determines all the
finite dimensional marginals of the distribution of
$\nu$.\footnote{We recall a standard fact: if  \,$U_1,\dots, U_l$\, are
positive random variables such that
$U_1+\cdots+ U_l\in[0,1]$ almost surely then, for Lebesgue-almost all
\,$(a_1,\dots, a_l)\in\R^l$,
$$\P(U_1<a_1,\dots, U_l<a_l)
=\lim_{n\to\infty}
\sum_{{{k_0,\dots, k_l\geq 0}
\atop{{k_0+\cdots+ k_l=n}\atop{ {k_1\over n}<a_1,\dots, {k_l\over
n}<a_l}}}}{n!\over k_0!\cdots k_l!}\E[(1-U_1-\cdots-U_l)^{k_0}\cdot
U_1^{k_1}\cdots U_l^{k_l}].$$
}

If card\,$T=1$, we get all the moments involving the $\nu_r$'s. And
if $t$ is the
unique element of $T$, then $\nu_t=1-\sum_{r\in R}\nu_r$, and we get all the
moments of $\nu$.
\qed

When card\,$T\geq2$  the law of the environment can be determined
in some cases, but not in general. (Accordingly, Corollary 1 of
\cite{ES} should be amended; it holds in fact if card $T\leq1$,
but not in complete generality.) Here are two examples:

\noindent{\sl{E\tiny{XAMPLE} \ }\small1}.\, We consider the two
following walks on $\Z$:

\noindent$\bullet$\ \  The first one has a deterministic
environment, and moves from $x\in \Z$  with probability $1\over2$
to $x+1$,  with probability $1\over2$ to $x+2$.

\noindent$\bullet$\ \  The environment in the second walk is
coin-tossed independently at each site $x\in \Z$: with probability
${1\over2}$, the transition probability to $x+1$ is 1, and with
probability ${1\over2}$, the transition probability to $x+2$ is 1.

Here $T=\{1,2\}$, and, obviously, the  two walks have the same
law.

\noindent{\sl{E\tiny{XAMPLE} \  }\small 2}.\, Again, \,$G=\Z$; %
 for any
$x\in\Z$, with probability ${1\over 2}$, the transition
probability from $x$ to  $x$ is  equal to ${1\over2}$, and
the transition
probability from $x$ to  $x+1$ is  also equal to ${1\over2}$; %
and, with
probability ${1\over 2}$, the transition probability from $x$
to $x+2$ is
equal to 1. In this case, $T=\{1,2\}$, but the
distribution of the environment is almost surely completely determined
by the single trajectory we observe.

We can only see 0-transitions (from a site $x$ to itself), 1-transitions ($x\to x+1$) 
and 2-transitions ($x\to x+2$). So $\mu$, which is the law of $\nu(0)$, 
satisfies  $\mu(\nu_0+\nu_1+\nu_2=1)=1$. The fact that a 0-transition is never followed %
by a 2-transition tells us that if $\nu_2>0$, then $\nu_0=0$. Statistics on sites from which there 
are 0-transitions or 1-transitions reveals that the (conditional) %
distribution of the 
number of 0-transitions from such a site  
is geometric. But a geometric distribution cannot be a nondegenerate 
convex combination of geometric distributions, and the conditional 
number of 0-transitions from a visited site $x$ for which $\nu(x,0)$ 
is given (and is in $]0,1[$) has a geometric distribution. We deduce that 
there is exactly one value $a$ ($={1\over 2}$) such that almost surely, for all $x\in \Z$, 
if $\nu(x,0)\neq 0$, then $\nu(x,0)=a$. And a simple computation shows 
that the proportion of visited sites from which a 0-transition is possible but 
does not take place fits exactly with the proportion of visited sites 
from which the first (and unique) transition is a 1-transition;
so $\mu(\nu_1>0 {\mathrm{\  and\  }} \nu_2>0)=0$. 
We deduce that $\mu$ is a convex
combination of
$\mu^1$ and $\mu^2$, where $\mu^1(\nu_0=\nu_1={1\over 2})=1$
and
$\mu^2(\nu_2=1)=1$. And the coefficients (both equal${1\over 2}$) in this convex combination are 
easily determined.

\subsection{Sampling iid trajectories\textcolor{white}{.}}~%

Here we show that infinitely many independent trajectories drawn
under $\P^\mu$ can be generated by concatenating steps of the
single trajectory at our disposal.

   For simplicity, we describe the
construction of just two trajectories, $X^1$ and $X^2$.
This is enough, since we can do the following: once  $X^1$ and $X^2$
are constructed, leave $X^1$ as it is and extract two
trajectories, say $X^3$ and $X^4$, out of $X^2$ exactly the
way we extracted  $X^1$ and $X^2$ out of $X$,
then, out of $X^4$,  get $X^5$ and $X^6$, and so on;
and the family \,$(X^1, X^3, X^5,\dots)$ is exactly what we want.

All we do is construct $X^1$ (resp. $X^2$) using in the ``natural" way the
transitions $\Delta_i^\n$ with $i$ odd (resp. even).
More formally, $X^1$ and $X^2$ are defined as follows:
$$ X^1_0=0,\quad  X^1_{n+1}=X^1_n+\Delta_{2\tau^1(n)-1}^{\vec{N}^{1}(n)},$$
where $ \vec{N}^{1}$ is to the process $X^1$ what $\vec{N}$ is to the
process $X$, and $$\tau^1(n):=
{\rm card}\{j\in\{0,\cdots,n\}\,|\,\vec{N}^1(j)=\vec{N}^1(n)\};$$
and, similarly,
$$ X^2_0=0,\quad  X^2_{n+1}=X^2_n+\Delta_{2\tau^2(n)}^{\vec{N}^{2}(n)},$$
where $ \vec{N}^{2}$ is to the process $X^2$ what $\vec{N}$ is to the
process $X$, and $$\tau^2(n):=
{\rm card}\{j\in\{0,\cdots,n\}\,|\,\vec{N}^2(j)=\vec{N}^2(n)\}.$$
The validity of this construction is an immediate consequence of Proposition 4.

\section{Remarks}

\subsection{Infinitude of \,$S$.}
The ``infinitude assumption" (according to which the random set $S$ of
sites visited by $X$ is almost surely
an infinite set) is made in order to avoid discussing rather
trivial cases. (If $S$ is finite, then precise knowledge of the
environment at {\emph{some}} sites is almost surely available;
but, unless some specific conditions are imposed on $\mu$,
complete knowledge of $\mu$ is out of reach.)

Proceeding along the general lines of the proof of Proposition 3, one can show
that the random set $S$ is either almost surely infinite or almost surely
finite.

\subsection{Countability of \,$G$.}
If we abandon the countability assumption on $G$, the set
$E=\{x\in G\,|\, P(X_1=x)>0\}$ remains countable, and our
sampling procedure  works just as well. Consequently, $\P^\mu$\, can be
determined in non-pathological situations (and in particular  if \,$G$\,is the
real line).\footnote{%
Of course, if \,$G$\, is countable, the \,$\sigma$-field  we  use
(without explicitly saying so) is the set of all subsets of
\,$G$;\, and if \,$G$\, is the real line, we take the Borel
\,$\sigma$-field on the line. But problems may arise if a
\,$\sigma$-field on \,$G$\, is not specified in advance and there
is no ``natural" \,$\sigma$-field on \,$G$:\, the very notion of
the law of \,$X$\, is problematic (and, in fact, even the notion
of the law of \,$X_1$\, does not make much sense). But even if a
\,$\sigma$-field on \,$G$\, is ``given", we aren't through. What
we want is to be able to determine, on the basis of the
observation of one realization of a random sequence \,$(U_1,
U_2,\dots)$ of iid random variables taking values in \,$G$,  the
probability distribution of \,$U_1$. Now if the \,$\sigma$-field
is generated by some countable \,$\pi$-system of subsets of \,$G$,
things are all right. Otherwise, there is no general positive
result.} (This can  also be seen by introducing a new kind of
reinforcement function which, to a given unordered history at a
site, associates the probability that the next transition falls
into some measurable set.) If there is no countable set
\,$D\subset G$ \,such that \,$X_1\in D$ \,almost surely, then
\,$\mu$ \,cannot be determined (as one can  see after studying the
first example of section 4). (If there is some countable set
\,$D\subset G$ \,such that \,$X_1\in D$ \,almost surely then,
almost surely, {\em{each}} \,$X_n$ \,is in the subgroup generated
by \,$D$;\, and since this subgroup is countable, all we did is
adaptable in an obvious way.)

\subsection{Structure of \,$G$.}
The choice of dealing with RWRE on a {\emph{group}} captures, we
think, the essence of the matter. We could have restricted
ourselves to groups like $\Z^d$ (or some other subgroups of
$\R^d$) without a substantial gain in simplicity. A group structure is
suitable (though not absolutely indispensable) if the notion of
iid random environment is to make sense. We could have dealt with
RWRE on homogenous spaces, or on some trees, without gaining new
insight.

\subsection{Assumptions on the environment.}
The requirement that environment at sites are iid can be loosened
in various ways.

\noindent{\sl{E\tiny{XAMPLE}.\ }}  $G=\Z$; there are two laws for
the environment at sites, say $\mu_0$ and $\mu_1$; environments at
sites are independent; and $\nu(n)$ is governed by $\mu_0$ if $n$
is even, by $\mu_1$ if $n$ is odd.

\noindent{\sl{E\tiny{XAMPLE}.\ }}  $G=\Z$; there are laws $\mu_0$,
$\mu_1$\dots for the environment at a site; $K$ is a random
variable taking values in the set $\{2, 3, \dots\}$; and,
conditioned on $K$, the $\nu(n)$ are independent and, for all $n$,
$\nu(n)$ is governed by $\mu_{n (\mathrm{mod}\,K)}$.

\noindent{\sl{E\tiny{XAMPLE}.\ }}  $G=\Z$; the couples %
$(\nu(2n),\nu(2n+1))\ \,(n\in\Z)$
are iid, but $\nu(0)$ and $\nu(1)$ are {\emph{not}} independent.

\subsection{Other reinforcements.} Our results on the determination
of the law of the process and on sampling iid trajectories can be
extended to various other edge-oriented reinforced random walks that do not
correspond to a random environment. Whenever an appropriate analogue of
Proposition 3 is valid, things work quite well. (A sufficient condition is
strict positivity of the restriction of $\,V_r\,$ to $\,(\Z_+^R)_0\,$ for all
\,$r\in R$.)


\begin{thebibliography}{10}
\bibitem{BK} Benjamini, I.,  Kesten, H., Distinguishing sceneries by observing
the scenery along a random walk path, J. Anal. Math. 69 (1996), 97--135.
\bibitem{ES}  Enriquez, N., Sabot, C., Edge-oriented reinforced random
walks and RWRE, C.R. Acad. Sci. Paris I 335 (2002), 941--946.
\bibitem{LM} L\"owe, M.,  Matzinger, H., Scenery reconstruction in two
dimensions with many colors. Ann. Appl. Probab. 12 (2002), no. 4, 1322--1347.
\bibitem{K} Kalikow, S. A., Generalized random walk in a random
environment. Ann.
Probab. 9 (1981), no. 5, 753--768.
\bibitem{M} Matzinger, H.,  Reconstructing a three-color scenery by
observing it
along a simple random walk path. Random Structures Algorithms 15 (1999), no. 2,
196--207.

\end{thebibliography}
\end{document}